\documentclass[twoside, 11pt]{article}

\usepackage{mathrsfs,amsfonts,amsmath}
\RequirePackage{ifthen,calc}
\usepackage{theorem}
\usepackage[dvips]{color}

\def\R{{\mathbb R}}

\def\E{{\mathbb E}}

\def\N{{\mathbb N}}

 \catcode`@=11
 \def\@evenhead{\hbox to\textwidth{\footnotesize\rm\thepage \hfill
  {\it Dai H.}}} % authors name

 \def\@oddhead{\hbox to \textwidth{\footnotesize{\it Approximation to RL-MBS
 } \hfill\thepage}}% abbreviate title

 \renewcommand{\section}{\makeatletter
 \renewcommand{\@seccntformat}[1]{{\csname the##1\endcsname.}\hspace{0.45em}}
 \makeatother \@startsection
{section}%                                            the name
{1}%                                                  the level
{0pt}%                                                the indent
{\baselineskip}%                                      the beforeskip
{0.5\baselineskip}%                                   the afterskip
{\normalsize\bfseries\mathversion{bold}}}

\renewcommand{\subsection}{\makeatletter
 \renewcommand{\@seccntformat}[1]{{\csname the##1\endcsname.}\hspace{0.45em}}
 \makeatother \@startsection
{subsection}%                                            the name
{1}%                                                  the level
{0pt}%                                                the indent
{\baselineskip}%                                      the beforeskip
{0.5\baselineskip}%                                   the afterskip
{\normalsize\bfseries\mathversion{bold}}}
\catcode`@=12

\newtheorem{thm}{\noindent Theorem}[section]
\newtheorem{lem}{\noindent Lemma}[section]
\newtheorem{cor}{\noindent Corollary}[section]

\theoremheaderfont{\normalfont\bfseries}
\theorembodyfont{\slshape} {\theorembodyfont{\rmfamily}} {\theorembodyfont{\rmfamily}
\newtheorem{defn}{\noindent Definition}[section]}
{\theorembodyfont{\rmfamily} } {\theorembodyfont{\rmfamily}
}
{\theorembodyfont{\rmfamily} } {\theorembodyfont{\rmfamily}
}
{\theorembodyfont{\rmfamily} } {\theorembodyfont{\rmfamily} } {\theorembodyfont{\rmfamily}
}
 \def\beqlb{\begin{eqnarray}}\def\eeqlb{\end{eqnarray}}
 \def\beqnn{\begin{eqnarray*}}\def\eeqnn{\end{eqnarray*}}
 \numberwithin{equation}{section}
\def\qed{\hfill$\square$\smallskip}

\begin{document}
\title{  Approximation to  multifractional Riemann-Liouville Brownian sheet}
\author{\small Hongshuai Dai\thanks{E-mail:  mathdsh@gmail.com}
 \\ \small  School of Mathematics and Information Sciences, Guangxi University,
 \\ \small Nanning, 530004  China
}
\maketitle
\begin{abstract}
\noindent In this paper, we first introduce  multifrational Riemann-Liouville Brownian sheets. Then, we show a result of  approximation in law of the multifractional Riemann-Liouville Brownian sheet. The construction of these approximations is based on  a sequence of I.I.D random variables.
\end{abstract}



\section{Introduction}

A one-dimensional fractional Brownian  motion (FBM) $\xi_H=\{\xi_H(t);\;t\in\R_+\}$ with Hurst index $H\in(0,\;1)$ is a real-valued, centered Gaussian process with covariance function given by
\beqlb\label{ds1-1}
\E\Big[\xi_H(t)\xi_H(s)\Big]=\frac{1}{2}\Big[|s|^{2H}+|t|^{2H}-|t-s|^{2H}\Big],\;s,\;t\in\R_+.
\eeqlb
It was introduced, as a moving-average Gaussian process, by Mandelbrot and Van Ness \cite{MN1968}.

Fractional Brownian motion possesses interesting properties such as self-similarity of order $H\in(0,\;1)$, stationary increments and long-range dependence (when $H>\frac{1}{2}$),  which  make it a good candidate for modeling different phenomena in,  for example, finance and telecommunication. However, it does not represent a casual time-invariant system as there is no well-defined impulse response function. Barnes and Allan \cite{BA1966} introduced the  fractional  Riemann-Liouville (RL) Brownian motion (RL-FBM) $V_H=\{V_H(t);\; t\in\R_+\}$ based on the following Riemann-Liouville fractional integral,
\beqlb V_H(t)=\frac{1}{\Gamma(H+\frac{1}{2})}\int_{\R_+}
(t-u)_+^{H-\frac{1}{2}}dB(u),
\eeqlb where $B(u)$ is a standard Brownian motion, $x_+=\max\{x,0\}$, and $\Gamma$ is the gamma function.

$V_H$  represents a linear system driven by white noise
with the impulse response function $\frac{t^{H-\frac{1}{2}}}{\Gamma(H+\frac{1}{2})}$. $V_H$  shares with $\xi_{H}$ many properties which include self-similarity, regularity of sample paths, etc.- with one notable exception that its increment process   is non stationary. But it was shown by Lim \cite{L2001} that the increment process  satisfies some weaker forms of stationary. This provides some flexibility necessary for practical applications (see Lim \cite{L2001} for more information). However, this model may be restrictive, due to the fact that all of its regularity and fractal properties are governed by the single Hurst parameter $H$. In order to overcome this limitation,  Lim and Muniandy \cite{LM2002} extended  RL-FBM to  RL-multifractional Brownian motion (RL-MBM).
Roughly speaking, RL-MBM is governed by a Hurst function $H(t)$ with certain regularity in place of the constant Hurst parameter $H$ in RL-FBM. We refer  to Lim \cite{L2001} and Lim and Muniandy \cite{LM2000,LM2002} for properties of RL-MBMs.

There are two typical multiparameter extensions of FBMs. One is the L$\acute{e}$vy's fractional Brownian random field with parameter $H\in(0,\;1)$ (see Ciesielski and Kamont \cite{CK1995}) : a centered Gaussian process $Y$ with covariance function given by
\beqnn
\E[Y(t)Y(s)]=\frac{1}{2}\Big[||s||^{2H}+||t||^{2H}-||t-s||^{2H}\Big],\;s,t\in\R_+^d,
\eeqnn
 where $||\cdot||$ denotes the Euclidean norm.

The other multiparameter extension of the fractional Brownian motion is the anisotropic fractional Brownian sheet introduced by Kamont \cite{K1996}.  Since then, fractional Brownian sheets have been studied extensively as a representative of anisotropic Gaussian random fields. See, for example, Wu and Xiao \cite{W2007} and the references therein for further information. Still, the regularity of fractional Brownian sheets does not evolve in the $d-$ dimensional time parameter $t\in\R_+^d$.

To model the anisotropic Gaussian random fields whose regularity evolves in time, such as images, Ayache and L$\acute{e}$ger \cite{AL2000} and Herbin \cite{H2006} introduced so-called multifractional Brownian sheets (MFBS) in terms of their moving average representations and harmonizable representations, where the constant Hurst vector of fractional Brownian sheets is substituted by Hurst functionals. We refer to Ayahce and L$\acute{e}$ger \cite{AL2000}, Herbin \cite{H2006}, Meerschaert, Wu and Xiao \cite{M2008} and the references therein for more information.

Now, we are ready to introduce the  multifractional Riemann-Liouville Brownian sheet (RL-MBS).
\begin{defn}\label{ds1-defn2}
Let $H(t)=(H_1(t),\cdots,H_d(t))$ be a measurable function in $t\in\R^d_+$ with values in $(0,\;1)^d$. A real-valued RL-multifractional Brownian sheet $X=\{X(t);\,t\in\R^d_+\}$ with functional Hurst index $H(t)$ is defined as the following Riemann-Liouville (RL) fractional integral
\beqlb\label{ds1-3}
X(t)=\int_{\R^d_+}\prod_{i=1}^d(t_i-u_i)_+^{H_i(t)-\frac{1}{2}}W(du),
\eeqlb
where $W(u)=\{W(u);u\in\R^d_+\}$ is a $d$-parameter Wiener process.
\end{defn}

Some authors have studied weak convergence to multifractional Brownian motions. Dai and Li \cite{D2010} presented a weak limit theorem for the multifractioan Brownian motion based on a Poisson process. By using Donsker's theorem, Dai \cite{D2012} showed an approximation of the RL-multifractional Brownian motion in Besov spaces. On the other hand,  Weak limit theorems for fractional Brownian sheets have attracted significant interest in recent years. For example, Bardina et al. \cite{Bardina2003} gave a result of approximation in law, in the space of the continuous functions on $[0,\;1]^2$, of two-parameter fractional Brownian sheet.  Bardina and Florit  \cite{B2008} showed a result of approximation in law of the $d$-parameter fractional Brownian sheet in the space of continuous functions on $[0,\;1]^d$.

Let us recall some known facts. Let $\big\{Z_k;\;k=(k_1,k_2,\dots,k_d)\in\N^d\big\}$ be an independent family of centered identically distributed random variables with variance 1.  Wichura \cite{Wichura1969} generalized Donsker's theorem and proved that the process
\beqlb\label{df1-4}
n^{\frac{d}{2}}\int_0^{t_d}\cdots\int_0^{t_1}\sum_{k\in\N^d}Z_k I_{[k-1,\;k)}(u\cdot n)du,
\eeqlb
where $I_{[k-1,\;k)}(u\cdot n)=I_{[k_1-1,\;k_1)\times\cdots\times [k_d-1,\;k_d)}(u_1n,\cdots,u_dn)$, converges in law to a $d$-parameter Wiener process.

Inspired by Wichura \cite{Wichura1969} and Bardina and Florit\cite{B2008}, it is natural to try to  approximate $X$ in law  by
\beqlb\label{ds1-5}
X_n(t)=\int_{\R^d_+}\Big[\prod_{i=1}^d(t_i-u_i)_+^{H_i(t)-\frac{1}{2}}\Big]\theta_n(u)du,
\eeqlb
where $\theta_n(u)=n^{\frac{d}{2}}\sum_{k\in\N^d}Z_k I_{[k-1,\;k)}(u\cdot n)$.

In the rest of this paper, we assume  that $H(t)=\big(H_1(t),\cdots,H_d(t)\big)$ satisfies the following conditions:
\begin{itemize}
\item For every $i\in\{1,\cdots,d\}$, $0<\alpha_i \leq H_i(t)\leq \beta_i<1$ for all $t\in \R_+^d$.
\item $H(t)$ is $\gamma$-H$\ddot{o}$lder continuous, i.e.,
\beqnn
\big\|H(t)-H(s)\big\|\leq K \|t-s\|^\gamma, \;t,s\in \R_+^d,
\eeqnn
where $K>0$ is a constant.
\end{itemize}
 Furthermore, we assume that  $\E|Z_k|^m<\infty$  for any $k\in \N^d$ and $m\in \N$.

In this paper, we consider $[0,\;1]^d\subset \R^d_+$ with the usual partial order.
We will prove that the sequence of laws in $\mathcal{C}([0,\;1]^d)$ of the processes $\{X_n(t);t\in[0,\;1]^d\}$ defined by (\ref{ds1-5}) converges weakly to the law of $\{X(t);\; t\in [0,\;1]^d\}$ given by (\ref{ds1-3}).

Most of the estimates of this paper contain unspecified constants.
An unspecified positive and finite constant will be denoted by $K$,
which may not be the same in each occurrence. Sometimes we shall
emphasize the dependence of these constants upon parameters.

We end this section with some properties of RL-multifractional Brownian sheets.
\begin{lem}\label{ds1-lem1}
For any $s<t\in[0,\;1]^d $, there exits a constant $K>0$ such that
\beqnn
\E\Big[X(t)-X(s)\Big]^2 \leq K \|t-s\|^{2H},
\eeqnn
where $H=\min_{i\in\{1,\cdots,d\}}\{\alpha_i,\gamma\}$.
\end{lem}
 {\it Proof:}  First, we have
 \beqlb\label{ds1-a1}
 \E\Big[X(t)X(s)\Big]=\Pi_{i=1}^d\Big[\int_{[0,\;1]}(t_i-u_i)_+^{H_i(t)-\frac{1}{2}}(s_i-u_i)_+^{H_i(s)-\frac{1}{2}}du_i\Big].
 \eeqlb

 By (\ref{ds1-a1}),  we can consider the process $X=\{X(t)\}$ defined from $d$ independent RL-multifractional Brownian motions $Y^{(i)}$ with parameter $\tilde{H}_i(t_i)=H_i(t)$ by
 \beqlb\label{ds1-a2}
 Y(t)=Y^{(1)}(t_1)\cdots Y^{(d)}(t_d).
 \eeqlb

 It is easy to see that it has the same covariance function as the RL-mutifractional Brownian sheet $X$.

 From above arguments, we get
 \beqlb\label{ds1-a3}
 \E\Big[X(t)-X(s)\Big]^2&&=\E\Big[\Pi_{i=1}^dY^{(i)}(t_i)-\Pi_{i=1}^dY^{(i)}(s_i)\Big]^2\nonumber
 \\&&=\E\bigg[\Big(\Pi_{i=1}^dY^{(i)}(t_i)-Y^{(1)}(s_i)\Pi_{i=2}^dY^{(i)}(t_i)\Big)+\cdots \nonumber \\&&\qquad+\Big(\Pi_{i=1}^{d-1}Y^{(i)}(s_i)Y^{(d)}(t_d)-\Pi_{i=1}^dY^{(i)}(s_i)\Big)\bigg]^2.
  \eeqlb
Since $Y^{(i)}, i=1,\cdots,d$ are mutually independent  and $\E[Y^{(i)}]^2<K$,  we get that
\beqlb\label{ds1-a5}
\E\Big[X(t)-X(s)\Big]^2\leq K\sum_{i=1}^d\E\Big[Y^{(i)}(t_i)-Y^{(i)}(s_i)\Big]^2.
\eeqlb

 It follows from Dai \cite{D2012} that
 \beqlb\label{ds1-a6}
 \E\Big[Y^{(i)}(t_i)-Y^{(i)}(s_i)\Big]^2\leq K|t_i-s_i|^{2H_i(t)}+K|H_i(t)-H_i(s)|^2.
 \eeqlb

  Since $t,s\in[0,\;1]^d$,
 \beqlb\label{ds1-a7}
 |t_i-s_i|\leq 1.
 \eeqlb

 On the other hand, we have
 \beqlb\label{ds1-a8}
| H_i(t)-H_i(s)|<1.
 \eeqlb

By (\ref{ds1-a4}) to  (\ref{ds1-a8}), one can easily get that the lemma holds.\qed

It follows from the Kolmogorov continuity theorem  (see Theorem 2.3.1 in Khoshnevisan \cite{Khoshnevisan2002}) and Lemma \ref{ds1-lem1}  that
\begin{cor}\label{C-1}
$\{X(t)\}$ has a continuous modification.
\end{cor}

 The proof of Corollary \ref{C-1} is classical, so we omit the proof.

\section{Main result}
The main result of this paper is the following theorem.
\begin{thm}\label{ds1-thm1}
The sequence of laws in $\mathcal{C}([0,\;1]^d)$ of the processes $\{X_n(t);t\in[0,\;1]^d\}$ defined by (\ref{ds1-5}) converges weakly to the law of $\{X(t);\; t\in [0,\;1]^d\}$ given by (\ref{ds1-3}).
\end{thm}

In order to prove Theorem \ref{ds1-thm1}, we have to check that the family of laws of the processes $\{X_n\}$ is tight. First, we need  a technical lemma.

\begin{lem}\label{ds2-lem1} For any even $m\in\N$, there exists a constant $K>0$ such that for any nonnegative functions $f_i\in L^2([0,\;1])$, $i=1,\cdots,d$,
\beqlb\label{ds2-1}
\E\Bigg[\int_{[0,\;1]^d}\Big[\prod_{i=1}^d f_i(u_i)\Big]\theta_n(u)du\Bigg]^m\leq K\prod_{i=1}^d \Big[\int_0^1 f^2_i(u_i)du_i\Big]^{\frac{m}{2}}.
\eeqlb
\end{lem}

{\it Proof:}  We have
\beqlb\label{ds2-2}
\E\Big[\int_{[0,\;1]^d}\big[\prod_{i=1}^d f_i(u_i)\big]\theta_n(u)du\Big]^m&&=\E\bigg[\int_{[0,\;1]^{dm}}\big[\prod_{j=1}^m\prod_{i=1}^d f_i(u_i^j)\big]\big[\prod_{j=1}^m\theta_n(u^j)\big]du^1\cdots du^m\bigg]\nonumber
\\&&=\int_{[0,\;1]^{dm}}\Big[\prod_{j=1}^m\prod_{i=1}^d f_i(u_i^j)\Big]\E\Big[\prod_{j=1}^m\theta_n(u^j)\Big]du^1\cdots du^m.
\eeqlb
 On the other hand, we have
\beqlb\label{ds2-3}
\E\Big[\prod_{j=1}^m\theta_n(u^j)\Big] &&=n^{\frac{dm}{2}}\E\Big[\prod_{j=1}^m \Big(\sum_{k\in \N^d}Z_{k}I_{[k-1,\;k)}(u^j\cdot n)\Big)\Big]\nonumber
\\&&= n^{\frac{dm}{2}}\E\Big[\sum_{k^1,\cdots,k^m\in \N^{d}}Z_{k^1}\cdots Z_{k^m}I_{[k^1-1,\;k^1)}(u^1\cdot n)\cdots I_{[k^m-1,\;k^m)}(u^m\cdot n)\Big]\nonumber
\\&&\leq K  n^{\frac{dm}{2}}\sum_{(k^1,\cdots,k^m)\in A_m}I_{[k^1-1,\;k^1)}(u^1\cdot n)\cdots I_{[k^m-1,\;k^m)}(u^m\cdot n),
\eeqlb
where
\beqnn
A_m&&=\Big\{(k^1,\cdots,k^m)\in \N^{dm}; \;\textrm{for all}\;i\in\{1,\cdots,m\},
\\ &&\quad\qquad \textrm{there exists}\; j\in\{1,\cdots,m\}\setminus\{i\} \textrm{such that}\; k^i=k^j\Big\},
\eeqnn
since $Z_k, k\in\N^d$ are independent and $\E|Z_k|^m<\infty$ for all $m\in\N $ and $k\in \N^d$.

We have
\beqlb\label{ds2-4}
\sum_{(k^1,\cdots,k^m)\in A_m}I_{[k^1-1,\;k^1)}(u^1\cdot n)\cdots I_{[k^m-1,\;k^m)}(u^m\cdot n)\leq I_{D_m}(u^1,\cdots,u^m),
\eeqlb

where
\beqnn
D_m&&=\Big\{(u^1,\cdots,u^m)\in [0,\;1]^{dm}; \;\textrm{for all}\;i\in\{1,\cdots,m\},
\\ &&\quad\qquad \textrm{there exists}\; j\in\{1,\cdots,m\}\setminus\{i\}\; \textrm{such that for all}\; k=1,\cdots,d,
\;|u^i_k-u^j_k|\leq \frac{1}{n}
\Big\}.
\eeqnn
In fact, we can bound $I_{D_m}$ by a finite sum of  products of indictators. Moreover all the $d m$ variables, $u^1_1,\cdots,u^m_d$, appear in each product of indicators, but each variable  appears only in each one of indicators of each product. Then we can bound (\ref{ds2-2}) by a sum of products of the following two kinds of terms:
\begin{itemize}
\item For any $i\neq j$, $i,j\in\{1,\cdots,m\}$,
\beqlb\label{ds2-5}
\tilde{I}_1=Kn^d\int_{[0,1]^{2d}}\prod_{k=1}^d \big[f_k(u_k^i)f_k(u_k^j) \big] \prod_{k=1}^d I_{[0,\;\frac{1}{n})}(|u_k^i-u_k^j|) du^idu^j.
\eeqlb
\item For any $i\neq j\neq r$, $i,j,r\in\{1,\cdots,m\}$,
\beqlb\label{ds2-6}
&&\tilde{I}_2=Kn^{\frac{3d}{2}}\int_{[0,1]^{3d}}\prod_{k=1}^d \big[f_k(u_k^i)f_k(u_k^j) f(u_k^r)\big]\nonumber
\\&& \qquad\qquad\prod_{k=1}^d\Big[I_{[0,\;\frac{1}{n})}(|u_k^i-u_k^j|)I_{[0,\;\frac{1}{n})}(|u_k^i-u_k^r|)I_{[0,\;\frac{1}{n})}(|u_k^j-u_k^r|)\Big]du^idu^rdu^j.\quad
\eeqlb
\end{itemize}

In order to prove the lemma, we only need to show
\beqlb\label{ds2-7}
\tilde{I}_1\leq K\int_{[0,\;1]^d}\prod_{k=1}^d\big[{f^2_k(u_k)}\big]du,
\eeqlb
and
\beqlb\label{ds2-8}
\tilde{I}_2\leq K \Big[ \int_{[0,\;1]^d}\prod_{k=1}^d\big[{f^2_k(u_k)}\big]du\Big]^{\frac{3}{2}}.
\eeqlb
First, we show that (\ref{ds2-7}) holds. $\tilde{I}_1$ can be bounded by
\beqlb\label{ds2-9}
 &&Kn^{d}\int_{[0,\;1]^{2d}} \prod_{k=1}^d \big[f_k(u^j_k)\big]^2 \prod_{k=1}^d I_{[0,\;\frac{1}{n})}(|u^j_k-u^i_k|)du^jdu^i\nonumber
\\&&\quad\qquad+K n^{d}\int_{[0,\;1]^{2d}} \prod_{k=1}^d \big[ f_k(u^i_k)\big]^2 \prod_{k=1}^d I_{[0,\;\frac{1}{n})}(|u^j_k-u^i_k|)du^jdu^i.
\eeqlb
where we have used the elementary inequality $2ab\leq a^2+b^2$.

The two terms in (\ref{ds2-9}) can be done with the same method. So we only deal with the first one.
In fact,
\beqlb\label{ds2-10}
&&Kn^{d} \int_{[0,\;1]^{2d}} \prod_{k=1}^d \big[f_k(u^j_k)\big]^2 \prod_{k=1}^d I_{[0,\;\frac{1}{n})}(|u^j_k-u^i_k|)du^idu^j\nonumber
\\&& \qquad\qquad \leq K \int_{[0,\;1]^{d}}\prod_{k=1}^d \big[f_k(u^j_k)\big]^2 du^j.
\eeqlb
From (\ref{ds2-9}) and (\ref{ds2-10}), we can get that (\ref{ds2-7}) holds.

Now, we prove (\ref{ds2-8}).  $\tilde{I}_2$ equals
\beqnn\label{ds2-12}
&&Kn^{\frac{3d}{2}}\int_{[0,\;1]^{2d}}\Big[\prod_{k=1}^d \big[f_k(u_k^i)
I_{[0,\;\frac{1}{n})}(|u_k^i-u_k^j|)\big]\Big]\nonumber
\\&& \qquad\qquad\Big[\int_{[0,\;1]^{d}}\prod_{k=1}^d \big[f_k(u_k^r)f_k(u_k^j)I_{[0,\;\frac{1}{n})}(|u_k^i-u_k^r|)I_{[0,\;\frac{1}{n})}(|u_k^j-u_k^r|)\big]du^r\Big]du^idu^j.
\eeqnn
By  the H$\ddot{o}$lder inequality,
 \beqlb\label{ds2-11}
 \tilde{I}_2\leq Kn^{\frac{3d}{2}} G_1^{\frac{1}{2}}\times G_2^{\frac{1}{2}},
 \eeqlb  where
\beqnn
G_1=\int_{[0,\;1]^{2d}}\prod_{k=1}^d \Big[f_k(u_k^i)
I_{[0,\;\frac{1}{n})}(|u_k^i-u_k^j|)\Big]^2 du^idu^j,
\eeqnn
and
\beqnn
G_2&&=\int_{[0,\;1]^{2d}}\bigg[\int_{[0,\;1]^{d}}\prod_{k=1}^d \Big[f_k(u_k^r)f_k(u_k^j)I_{[0,\;\frac{1}{n})}(|u_k^i-u_k^r|)I_{[0,\;\frac{1}{n})}(|u_k^j-u_k^r|)\Big]du^r\bigg]^{2}du^idu^j.
\eeqnn

Using the same method as the proof of (\ref{ds2-7}), we can get
\beqlb\label{ds2-13}
G_1\leq K n^{-d}\int_{[0,\;1]^{d}}\prod_{k=1}^d f_k^2(u_k)du.
\eeqlb

Now, we deal with $G_2$.  $G_2$ equals

\beqlb\label{ds2-14}
&&\int_{[0,\;1]^{4d}}\prod_{k=1}^d \Big[f^2_k(u^j_k)f_k(u^{r_1}_k)f_k(u^{r_2}_k)I_{[0,\;\frac{1}{n})}(|u^{r_1}_k-u^j_k|)I_{[0,\;\frac{1}{n})}(|u^{r_2}_k-u^j_k|)\nonumber
\\&& \qquad\qquad I_{[0,\;\frac{1}{n})}(|u^{r_1}_k-u^i_k|)I_{[0,\;\frac{1}{n})}(|u^{r_2}_k-u^i_k|)\Big]du^{r_1}du^{r_2}du^{j}du^i.
\eeqlb

So, $G_2$ can be bounded by
\beqlb\label{ds2-15}
&&K n^{-d}\int_{[0,\;1]^{3d}}\prod_{k=1}^d \Big[f^2_k(u^j_k)f_k(u^{r_1}_k)f_k(u^{r_2}_k)I_{[0,\;\frac{1}{n})}(|u^{r_1}_k-u^j_k|)I_{[0,\;\frac{1}{n})}(|u^{r_2}_k-u^j_k|)\Big]du^{r_1}du^{r_2}du^{j}\nonumber
\\&&\quad\qquad\qquad\leq Kn^{-d}(G_3+G_4),
\eeqlb
where
\beqnn
G_3=\int_{[0,\;1]^{3d}}\prod_{k=1}^d \Big[f^2_k(u^j_k)f^2_k(u^{r_1}_k)I_{[0,\;\frac{1}{n})}(|u^{r_2}_k-u^j_k|)\Big]du^{r_1}du^{r_2}du^{j},
\eeqnn
and
\beqnn
G_4=\int_{[0,\;1]^{3d}}\prod_{k=1}^d \Big[f^2_k(u^j_k)f^2_k(u^{r_2}_k)I_{[0,\;\frac{1}{n})}(|u^{r_1}_k-u^j_k|)\Big]du^{r_1}du^{r_2}du^{j}.
\eeqnn
$G_3$ can be bounded by
\beqlb\label{ds2-16}
&&\int_{[0,\;1]^{2d}}\Big[\prod_{k=1}^d f^2_k(u^j_k)f^2_k(u^{r_1}_k)\Big]\Big[\int_{[0,\;1]^{d}}\prod_{k=1}^d I_{[0,\;\frac{1}{n})}(|u^{r_2}_k-u^j_k|)du^{r_2}\Big]du^{r_1}du^{j}\nonumber
\\&& \quad\qquad \leq Kn^{-d}\int_{[0,\;1]^{2d}}\prod_{k=1}^d \big[f^2_k(u^j_k)f^2_k(u^{r_1}_k)\big]du^{r_1}du^{j}\nonumber
\\&&\quad\qquad =Kn^{-d}\Big[\int_{[0,\;1]^{d}}\prod_{k=1}^d \big[f^2_k(u_k)\big]du\Big]^2.
\eeqlb

Using the same method as above, we can get that $G_4$ can be bounded by
\beqlb\label{ds2-17}
Kn^{-d}\bigg[\int_{[0,\;1]^{d}}\prod_{k=1}^d\Big[ f^2_k(u_k)\Big]du\bigg]^2.
\eeqlb

It follows from (\ref{ds2-12}) to (\ref{ds2-17}) that the inequality (\ref{ds2-8}) holds. This completes the proof. \qed

\begin{lem}\label{ds2-lem2}
For any even $m\in\N$, there exists a constant $K>0$ such that for any $s,t\in[0,\;1]^d$ with $s<t$,
\beqnn
\E\Big[X_n(t)-X_n(s)\Big]^m \leq K\big\| t-s\big\|^{mH},
\eeqnn
where $H=\min_{i\in\{1,\cdots,d\}}\{\alpha_i,\;\gamma\}$.
\end{lem}

{\it Proof:} It follows from Lemma \ref{ds2-lem1} that
\beqnn
\E\Big[X_n(t)-X_n(s)\Big]^m &&\leq K \bigg[\int_{[0,\;1]^d}\prod_{i=1}^d \Big((t_i-u_i)_+^{H_i(t)-\frac{1}{2}}-(s_i-u_i)_+^{H_i(s)-\frac{1}{2}}\Big)^2du\bigg]^{\frac{m}{2}}
\\&&\leq K \| t-s\|^{mH}.
\eeqnn
\qed

On the other hand,  one can easily prove that
\begin{cor}
The process $X_n$ is continuous for every $t\in [0,\;1]^d$  a.s.
\end{cor}

Next, we will prove the main result of this paper.

{\it Proof of Theorem \ref{ds1-thm1}:} By Bickel and Wichura \cite{Bichel1997}, in order to prove  Theorem \ref{ds1-thm1}, we only need to  prove that the family of laws of processes $\{X_n(t)\}$ is tight, and identify the limit law of any convergent subsequence as the law of the $d$-parameter RL-multifractional Brownian sheet $X$. The tightness comes from Lemma \ref{ds2-lem2} and Bickel and Wichura \cite{Bichel1997}.

Now, we proceed to identify the limit law. It is sufficient to show that for any $q\in\N, a_1,\cdots,a_q\in\R$ and $t^1,\cdots,t^q\in[0,\;1]^d$,
\beqlb\label{ds2-18}
\E\bigg[\exp\Big\{i\xi \sum_{j=1}^qa_jX_n(t^j)\Big\}\bigg]\to\E\bigg[\exp\Big\{i\xi \sum_{j=1}^qa_jX(t^j)\Big\}\bigg],
\eeqlb
as $n\to\infty$.

For any $j\in\{1,\cdots,q\}$ and $i\in\{1,\cdots,d\}$, we have
\beqlb\label{R-1}
(t_i^j-u_i)_+^{H_i(t^j)-\frac{1}{2}}\in L^2([0,\;1]),
\eeqlb
and
\beqlb\label{R-2}
(t_i^j-u_i)_+^{H_i(t^j)-\frac{1}{2}}\geq 0.
\eeqlb

Therefore, for any $j\in\{1,\cdots,q\}$ and $i\in\{1,\cdots,d\}$, we can find a sequence $\{\rho^{j}_{i,k}(u_i)\}_{k\in\N}$ of simple functions such that
\beqlb\label{ds2-20}
\int_{[0,\;1]}\Big(\rho^j_{i,k}(u_i)-(t_i^j-u_i)_+^{H_i(t^j)-\frac{1}{2}}\Big)^2du_i\to 0,\;\textrm{as}\;k\to\infty,
\eeqlb

and \beqlb\label{R-3}(t_i^j-u_i)_+^{H_i(t^j)-\frac{1}{2}}-\rho^{j}_{i,k}(u_i)\geq 0,\;\textrm{for all}\; k\in\N.
\eeqlb

In order to simplify the notation, we define
\beqnn
X_n^{j,k}=\int_{[0,\;1]^d} \big[\prod_{i=1}^d\rho_{i,k}^j(u_i)\big]\theta_n(u)du,
\eeqnn
and
\beqnn
X^{j,k}=\int_{[0,\;1]^d} \prod_{i=1}^d\rho_{i,k}^j(u_i)W(du).
\eeqnn

Note that

\beqlb\label{1}
\Bigg |\E\bigg[\exp\Big\{i\xi \sum_{j=1}^qa_jX_n(t^j)\Big\}\bigg]-\E\bigg[\exp\Big\{i\xi \sum_{j=1}^qa_jX(t^j)\Big\}\bigg]\Bigg|\leq I_1+I_2+I_3,
\eeqlb
where
\beqnn
I_1&&=\bigg|\E\Big[\exp\Big\{i\xi \sum_{j=1}^qa_jX_n(t^j)\Big\}-\exp\Big\{i\xi \sum_{j=1}^qa_jX_n^{j,k}\Big\}\Big]\bigg|,
\\I_2&&=\bigg|\E\Big[\exp\Big\{i\xi \sum_{j=1}^qa_jX(t^j)\Big\}-\exp\Big\{i\xi \sum_{j=1}^qa_jX^{j,k}\Big\}\Big]\bigg|,
\eeqnn
and
\beqnn
I_3&&=\bigg|\E\Big[\exp\Big\{i\xi \sum_{j=1}^qa_jX_n^{j,k}\Big\}-\exp\Big\{i\xi \sum_{j=1}^qa_jX^{j,k}\Big\}\Big]\bigg|.
\eeqnn

We first deal with $I_1$. By the mean value theorem,
\beqlb\label{ds2-19}
I_1\leq K\max_{1\leq j\leq q} \Big\{\E\Big|X_n(t^j)-X_n^{j,k}\Big|\Big\}
\eeqlb
On the other hand, for $j\in\{1,2,\cdots,q\}$
\beqlb\label{ds2-21}
&&\E\Big|X_n(t^j)-X_n^{j,k}\Big|\nonumber
\\&&\qquad= \E\Bigg[\bigg|\int_{[0,\;1]^d}\Big[(t^j_1-u_1)_+^{H_1(t^j)-\frac{1}{2}}-\rho_{1,k}^j(u_1)\Big]\Big[\prod_{i=2}^d(t_i^j-u_i)_+^{H_i(t^j)-\frac{1}{2}}\Big]\theta_n(u)dudu\nonumber
\\&&\qquad+\int_{[0,\;1]^d}\rho_{1,k}^j(u_1)\Big[\prod_{i=2}^d(t_i^j-u_i)_+^{H_i(t^j)-\frac{1}{2}} -\prod_{i=2}^d \rho_{i,k}^j (u_i)\Big]\theta_n(u)du\bigg|\Bigg]\nonumber
\\&&\qquad\leq \E\Bigg[\bigg|\int_{[0,\;1]^d}\Big[(t^j_1-u_1)_+^{H_1(t^j)-\frac{1}{2}}-\rho_{1,k}^j(u_1)\Big]\Big[\prod_{i=2}^d(t_i^j-u_i)_+^{H_i(t^j)-\frac{1}{2}}\Big]\theta_n(u)du\bigg|\Bigg]\nonumber
\\&&\qquad+\E\Bigg[\bigg|\int_{[0,\;1]^d} \rho_{1,k}^j(u_1)\Big[\prod_{i=2}^d(t_i^j-u_i)_+^{H_i(t^j)-\frac{1}{2}} -\prod_{i=2}^d \rho_{i,k}^j(u_i) \Big]\theta_n(u)du\bigg|\Bigg].
\eeqlb

On the other hand,
\beqlb\label{ds2-22}
&&\E\Bigg[\bigg|\int_{[0,\;1]^d}\Big[(t^j_1-u_1)_+^{H_1(t^j)-\frac{1}{2}}-\rho_{1,k}^j(u_1)\Big]\Big[\prod_{i=2}^d(t_i^j-u_i)_+^{H_i(t^j)-\frac{1}{2}}\Big]\theta_n(u)du\bigg|\Bigg]\nonumber
\\&&\quad\qquad\leq K \Bigg[\int_{[0,\;1]^d}\bigg[\Big((t^j_1-u_1)_+^{H_1(t^j)-\frac{1}{2}}-\rho_{1,k}^j(u_1)\Big)\prod_{i=2}^d(t_i^j-u_i)_+^{H_i(t^j)-\frac{1}{2}}\bigg]^2du\Bigg]^{\frac{1}{2}}\nonumber
\\&&\quad\qquad \leq K \bigg[\int_{[0,\;1]}\Big((t^j_1-u_1)_+^{H_1(t^j)-\frac{1}{2}}-\rho_{1,k}^j(u_1)\Big)^2du_1\bigg]^{\frac{1}{2}},
\eeqlb
where we have used the H$\ddot{o}$lder inequality and Lemma \ref{ds2-lem1}.

It follows from (\ref{ds2-20}) and (\ref{ds2-22}) that when $k\to\infty$,
\beqlb\label{ds2-23}
\E\Bigg[\bigg|\int_{[0,\;1]^d}\Big[(t^j_1-u_1)_+^{H_1(t^j)-\frac{1}{2}}-\rho_{1,k}^j(u_1)\Big]\Big[\prod_{i=2}^d(t_i^j-u_i)_+^{H_i(t^j)-\frac{1}{2}}\Big]\theta_n(u)du\bigg|\Bigg]\stackrel{U}{\to} 0,
\eeqlb
where $\stackrel{U}{\to}$ denotes uniform convergence with respect to $n$.

Repeating the above procedure,  we can get
$$
\E\Big|X_n(t^j)-X_n^{j,k}\Big|\stackrel{U}{\to}0 , \;\textrm{as}\; k\to\infty.
$$
Therefore,
\beqlb\label{ds2-24}
I_1\stackrel{U}{\to}0, \;\textrm{as}\; k\to\infty.
\eeqlb

Similarly, by the mean value theorem, we can get
\beqlb\label{R-4}
I_2\leq K \max_{1\leq j\leq q}\Big\{\E\big|X(t^j)-X^{j,k}\big|\Big\}
\eeqlb
Applying the H$\ddot{o}$lder inequality and the properties of the stochastic integral, we can get
\beqlb\label{ds2-25}
\E\big|X(t^j)-X^{j,k}\big|\leq K\Bigg[\int_{[0,\;1]^d}\Big(\prod_{i=1}^d\rho^{j}_{i,k}(u_i)-\prod_{i=1}^d(t_i^j-u_i)_+^{H_i(t^j)-\frac{1}{2}}\Big)^2du\Bigg]^{\frac{1}{2}}.
\eeqlb
Using the same method as the proof of (\ref{ds2-24}), we can get that
\beqlb\label{ds2-26}
I_2\stackrel{U}{\to}0, \;\textrm{as}\; k\to\infty.
\eeqlb

Due to the result established by Wichura \cite{Wichura1969}, the law of  the process
\beqnn
\int_{0}^{t_1}\cdots\int_{0}^{t_d}\theta_n(u)du
\eeqnn
converges weakly to the law of  a $d$-parameter Wiener process.  Then,  the linear combinations of increments of
\beqnn
\int_{0}^{t_1}\cdots\int_{0}^{t_d}\theta_n(u)du
\eeqnn
will converge in law to the same combinations of increments of the $d$-parameter Wiener process. So
\beqlb\label{ds2-27}
I_3\to 0,
\eeqlb
as $n\to \infty$.

Combining (\ref{ds2-18}), (\ref{1}), (\ref{ds2-24}), (\ref{ds2-26}) and (\ref{ds2-27}), we can get that Theorem \ref{ds1-thm1} holds.\qed

\noindent{\bf Acknowledgment} The author would like to thank  an anonymous referee for the insightful suggestions and helpful comments. This work was supported by Natural Science Foundation of Guangxi, China (2012GXNSFBA053010).

\end{document}